\documentclass[12pt]{article}
\def\al{\alpha} \def\be{\beta} \def\ov{\over} \def\si{\sigma}
\def\sp{\vskip2ex} \def\Ga{\Gamma} 
\def\ph{\varphi}

\begin{document}

\begin{center}{\bf Two Elementary Derivations of the Pure Fisher-Hartwig Determinant}\end{center}
\begin{center}{\bf Albrecht B\"ottcher and Harold Widom
}\end{center}

\noindent
By the ``pure Fisher-Hartwig determinant'' we mean the Toeplitz determinant
$D_n(\ph):={\rm det} (\ph_{i-j})_{i,j=1}^n$ where
the $\ph_k$ are the Fourier coefficients of
\[\ph (z)=(1-z)^\al\,(1-z^{-1})^\be,\]
a so-called pure Fisher-Hartwig singularity.
The $k$th Fourier coefficient of $\ph$ equals
\begin{equation}
(-1)^k{\Ga(\al+\be+1)\ov\Ga(\al+1-k)\,\Ga(\be+1+k)}.\label{coeff}
\end{equation}
The formula for the determinant
is
\begin{equation}
 D_n(\ph)=G(n+1)\ {G(\al+\be+n+1)\ov G(\al+\be+1)}\ {G(\al+1)\ov G(\al+n+1)}\
{G(\be+1)\ov G(\be+n+1)},\label{det}
\end{equation}
where $G$ is the Barnes $G$-function. This was deduced by Silbermann
and one of the authors
\cite{BS1} from a factorization of the
Toeplitz matrix $T_n(\ph)$ due to Duduchava and Roch.
Another proof was recently found by Basor and Chen \cite{BC} using the
theory of orthogonal polynomials, which motivated us to present the
two proofs of this note.

\vspace{3mm}
\noindent
{\bf First proof.} This proof
is analogous to the usual derivation of the Cauchy determinant
and its philosophy is that the most elegant way to determine a rational
function is to find its zeros and poles.

The factor $(-1)^k$ in (\ref{coeff}) will not affect the determinant.
We write the rest as
\[{\Ga(\al+\be+1)\ov\Ga(\al+1)\,\Ga(\be+1)}\;{\Ga(\al+1)\,\Ga(\be+1)\ov
\Ga(\al+1-k)\,\Ga(\be+1+k)}.\]
For the evaluation of $D_n(\ph)$ the first factor
will contribute in the end the factor
\begin{equation}
\left({\Ga(\al+\be+1)\ov\Ga(\al+1)\,\Ga(\be+1)}\right)^n.\label{factor}
\end{equation}
The remaining factor gives the determinant of the matrix $M$
with $i,j$ entry
\begin{equation}
 M_{ij}={\Ga(\al+1)\,\Ga(\be+1)\ov\Ga(\al+1-i+j)\,\Ga(\be+1+i-j)}.\label{M}
 \end{equation}

We think of $\det M$ as a function of $\al$, with $\be$ as
a parameter, and shall establish the following two facts:

\noindent(a)\ The only possible poles of $\det M$ (including $\infty$) are at
$-1,\ldots,-n+1$, with the pole at $-k$ having order at most $n-k$.

\noindent(b)\ For $k=1,\ldots,n-1$ $\det M$ has a zero at $\al=-\be-k$ of order
at least $n-k$.

Granting these for the moment, let us derive (\ref{det}). If $\det M$
had exactly the poles and zeros as stated it would be a constant depending on $\be$
times
\[\prod_{k=1}^{n-1}{(\al+\be+k)^{n-k}\ov(\al+k)^{n-k}}.\]
If there were more zeros or fewer poles, then in the representation of $\det M$
as a quotient of polynomials there would be at least one more non-constant factor in the
numerator than in the denominator. But then $\det M$ would not be analytic at
$\al=\infty$, which we
know it to be. Thus $\det M$ is a constant times the above. When $\al=0$ the matrix
is upper-triangular with diagonal entries all
equal to 1, so $\det M=1$ then. This determines the constant factor, and we deduce
\[\det M=\prod_{k=1}^{n-1}{k^{n-k}\,(\al+\be+k)^{n-k}\ov(\al+k)^{n-k}(\be+k)^{n-k}}.\]
Multiplying this by (\ref{factor}) gives (\ref{det}).
We now establish (a) and (b).

Proof of (a): The only possible finite poles of the $M_{ij}$ arise from
the poles of the numerator in (\ref{M}) at the negative integers $-k$. The pole at $-k$ will not be
cancelled by a pole in the denominator precisely when $j-i\ge k$. In particular for
there to be a pole we must have $k\le n-1$. The
order of the pole at $\al=-k$ in a term $\prod M_{i,\si(i)}$ in the expansion of
$\det M$ (here $\si$ is a permutation of $0,\ldots,n-1$) equals
\[\#\{i:\si(i)\ge i+k\}.\]
Since the inequality can only occur when $i<n-k$ the above number is
at most $n-k$. This establishes the statement about the possible finite poles. To see
that $\det M$ is analytic at $\al=\infty$, we observe that the order of the
pole of $M_{ij}$ there equals $i-j$. (The order is counted as negative when there is a zero.)
Hence the order of the pole there of $\prod M_{i,\si(i)}$ equals $\sum_i(i-\si(i))=0$.

Proof of (b): Let us write $M_{ij}(\al,\,\be)$ instead of $M_{ij}$ to
show its dependence on $\al$ and $\be$. A simple computation gives
for $i=1,\ldots,n-1$
\begin{eqnarray*}
& & M_{i,j}(\al,\,\be)+M_{i-1,j}(\al,\,\be)\\
& & =(\al+\be+1)\,{\Ga(\al+1)\,\Ga(\be+1)\ov\Ga(\al+2-i+j)\,\Ga(\be+1+i-j)}
={\al+\be+1\ov \al+1}\,M_{i,j}(\al+1,\be).
\end{eqnarray*}
In other words, if we add to each of the last $n-1$ rows of $M(\al,\,\be)$ the
preceding row we obtain $(\al+\be+1)/(\al+1)$ times the last $n-1$ rows of the
matrix $M(\al+1,\,\be)$. Then
we continue. If we apply these operations a total of $k$ times
the last $n-k$ rows of the matrix obtained from $M(\al,\,\be)$ in this way (which does not
change its rank) is equal to
\[{(\al+\be+1)\cdots(\al+\be+k)\ov (\al+1)\cdots(\al+k)}\]
times the last $n-k$ rows of the matrix $M(\al+k,\,\be)$. It follows that if we set
$\al=-\be-k$ in $M(\al,\,\be)$ we get a
matrix of rank at most $k$. From this it follows that if we differentiate $\det M(\al,\,\be)$
up to $n-k-1$ times with respect to $\al$ and set $\al=-\be-k$ we get zero. Thus there
is a zero there of order at least $n-k$.

\vspace{3mm}
\noindent
{\bf Second proof.} This proof does not aspire to elegance but is
rather the simple endeavor to go ahead straightforwardly.

Taking into account formula (\ref{coeff}) for the Fourier
coefficients of $\varphi$ we get
\[D_n(\varphi)=(\Gamma(\alpha+\beta+1))^n\,
{\rm det}\,\left(\frac{1}{\Gamma(\alpha+1-i+j)
\Gamma(\beta+1+i-j)}\right)_{i,j=1}^n.\]
Extracting the factor $1/\Gamma(\alpha+1+n-i)$ from the $i$th row
and $1/\Gamma(\beta+1+n-j)$ from the $j$th column,
we obtain
\begin{eqnarray}
& & \frac{D_n(\varphi)} {(\Gamma(\alpha+\beta+1))^n} =
\prod_{i=1}^n\frac{1}{\Gamma(\alpha+1+n-i)}\,
\prod_{j=1}^n\frac{1}{\Gamma(\beta+1+n-j)}\,D_n(\alpha,\beta)\nonumber\\[1ex]
& &
\hspace{31mm}=\frac{G(\alpha+1)}{G(\alpha+n+1)}\,
\frac{G(\beta+1)}{G(\beta+n+1)}\,
D_n(\alpha,\beta) \label{1}
\end{eqnarray}
with
\[D_n(\alpha, \beta)={\rm det}\left(
\prod _{\ell=1}^{n-j}(\alpha-i+j+\ell)\,
\prod_{k=1}^{n-i}(\beta+i-j+k)\,\right)_{i,j=1}^n.\]
The last row of $D_n(\alpha,\beta)$ is
\[\left(\quad\prod_{\ell=0}^{n-2}(\alpha-\ell)\quad
\prod_{\ell=0}^{n-3}(\alpha-\ell)\quad\ldots\quad
(\alpha-1)\alpha \quad\alpha\quad1\quad\right).\]
With the objective that the last row becomes
$(\; 0 \; 0 \; \ldots \; 0 \; 1\;)$,
we subtract $\alpha-n+2$ times column $2$ from column $1$,
$\alpha-n+3$ times column $3$ from column $2$, $\ldots\,$,
and finally $\alpha$ times column $n$ from column $n-1$.
What results is that
\[D_n(\alpha,\beta)=(n-1)!\,(\alpha+\beta+1)^{n-1}\,
D_{n-1}(\alpha+1, \beta).\]
Since $D_1(\alpha+n-1,\beta)=1$, it follows that

\begin{eqnarray}
& & D_n(\alpha,\beta) =  \prod_{k=1}^{n-1}(n-k)!\,
(\alpha+\beta+k)^{n-k}= \prod_{\ell=1}^n \Gamma(\ell)\,
\prod_{\ell=1}^n
\frac{\Gamma(\alpha+\beta+\ell)}{\Gamma(\alpha+\beta+1)}
\nonumber\\
& &
\hspace{16.5mm}= G(n+1)\,
\frac{G(\alpha+\beta+n+1)}{G(\alpha+\beta+1)}\,
\frac{1}{\big(\Gamma(\alpha+\beta+1))^n}.
\label{2}
\end{eqnarray}
Inserting (\ref{2}) in (\ref{1}) we arrive at the desired formula.

\pagebreak

\sp
\begin{center}{\bf Acknowledgment}\end{center}

\begin{center}
\begin{minipage}[t]{15cm}
The work of the second author was supported by National Science
Foundation grant DMS-0243982.
\end{minipage}
\end{center}

\sp\sp
\noindent
{\it Fakult\"at f\"ur Mathematik\\
TU Chemnitz\\
09107 Chemnitz, Germany\\
e-mail:} aboettch@mathematik.tu-chemnitz.de

\sp\sp
\noindent
{\it Department of Mathematics\\
University of California\\
Santa Cruz, CA 95064, USA\\
e-mail:} widom@math.ucsc.edu


\begin{thebibliography}{1}

\bibitem{BC} E. L. Basor and Y. Chen, {\it Toeplitz determinant from
compatibility conditions}, preprint.

\bibitem{BS1} A. B\"ottcher and B. Silbermann, {\it Toeplitz matrices and
determinants with Fisher-Hartwig symbols}, J. Funct. Anal. {\bf 62} (1985),
178--214.

\end{thebibliography}
\end{document}